\documentclass[12pt]{article}
\usepackage{amsmath}
\usepackage{amssymb}
\usepackage{amscd}
\usepackage{tabularx}
\usepackage{float}
\usepackage{theorem}
\usepackage{enumerate}
\def\eqnum#1{}

\begin{document}

\title
{Unilateral Small Deviations of Processes\\
Related to the Fractional Brownian Motion}

\date{}
\maketitle

\noindent {\footnotesize G.~MOLCHAN} \\

\noindent{\footnotesize
{\it International Institute of Earthquake Prediction Theory and Mathematical \\
Geophysics, Moscow}} \\
{\footnotesize
{\it The Abdus Salam International Centre
for Theoretical Physics, Trieste \\
}

\bigskip

\noindent {{\footnotesize \bf Abstract.}
{\footnotesize Let $x(s)$, $s\in R^d$ be a Gaussian
self-similar random process of index $H$. We consider the problem
of log-asymptotics for the probability $p_T$ that $x(s)$, $x(0)=0$
does not exceed a fixed level in a star-shaped expanding domain
$T\cdot \Delta$ as $T\to \infty$.
We solve the problem of the
existence of the limit,
$\theta :=\lim (-\log p_T)/(\log T)^D$, $T\to \infty$,
for the fractional Brownian sheet $x(s)$, $s\in [0,T]^2$ then $D=2$, and
we estimate $\theta$ for the integrated fractional Brownian motion
then $D=1$.}

\bigskip

\bigskip

\bigskip

\bigskip

\bigskip

\bigskip

\noindent \rule{7cm}{0.1mm}
\bigskip

\noindent {\footnotesize $^1$
Supported in part by the European Commission's Project
12975 (NEST) "Extreme Events: Causes and Consequences (E2-C2)".} \\

\bigskip

\noindent {\footnotesize {\it AMS 2000 Subject Classification}.
Primary - 60G15, 60G18}

\bigskip

\noindent {{\footnotesize \it Key words and phrases:}}
{\footnotesize small deviations, decay exponents,
fractional Brownian motion, level crossing probability}

\bigskip

\noindent {{\footnotesize \bf Abbreviated Title.}}
{\footnotesize Existence of Decay Exponents}
}}


\newpage

\section{Introduction}

Let $x(s)$, $s\in R^d$, $x(0)=0$ be a Gaussian random process
and $\Delta_T\ni 0$ an expanding subset of $R^d$. For many reasons it
is important for applications to know the asymptotics of the probability
$p_T$ that $x(s)$ does not exceed a fixed level in $\Delta_T$,

\begin{eqnarray}\label{qq1}
p_T = P\{ \sup x(s) < \mbox {\rm const},  s\in \Delta_T\}, \quad
T \to \infty
\eqnum{1}
\end{eqnarray}
(see [4], [6], [10, 11]). It has been sometimes possible to find the
order of logarithmic asymptotics of $p_T$ [4]. The second (in importance)
problem is the existence of the decay exponent, that is,

\begin{eqnarray}\label{qq2}
-\lim_{T\to \infty}\log p_T/\psi_T = \theta
\eqnum{2}
\end{eqnarray}
if $\log p_T=O(\psi_T)$. The exact values of $\theta$ are known only for
a few cases [12, 7].

Below we consider the fractional Brownian sheet $x(s)$, $s\in R^2$
for which

\begin{eqnarray}\label{qq3}
Ex(s)x(t) = B(s_1,t_1)B(s_2,t_2),
\eqnum{3}
\end{eqnarray}
where
$B(a,b)= \frac{1}{2}(|a|^{2H} + |b|^{2H} - |a-b|^{2H})$, $0<H<1$
is the correlation function of the fractional
Brownian motion $b_H(a)$, $a\in R^1$. In the case of the
Brownian sheet $(H=1/2)$ and $\Delta_T=[0,T]^2$ Cs\'aki et al. [1]
prove that

\begin{eqnarray*}
c_1/\rho_T \le -\log p_T/(\log T)^2 \le c_2 ,
\end{eqnarray*}
where $\rho_T=\log \log T$. Later, W.~Li and Q.~Shao [4] removed the
$\log \log$ term $\rho_T$ and formulated the problem to prove (2).
Below we do this for any $0<H<1$. With that end in view we consider the
problem (\ref{qq1},\ref{qq2}) for a
self-similar process of index $H$ ($H$-ss process). Self-similarity
means that
finite dimensional distributions of $\{ x(\lambda s)\}$ and
$\{ \lambda^H x(s)\}$ are identical for any $\lambda >0$.
For the $H$-ss process the problem (\ref{qq1}, \ref{qq2}) is reduced to a
similar one for events of the type
$\{ x(s)<0, s\in \Delta_T \setminus U_0\}$, where  $U_0$
is a vicinity of points $\{ s: Ex^2(s)=0\}$.
The existence of the decay exponent (\ref{qq2}) is deduced from
the positivity of (\ref{qq3}) and Slepian's lemma.
We consider also the integrated fractional
Brownian motion to find bounds of $\theta$ in this case.

\newpage

\section{General results}

In what follows $x(s)$, $s\in R^d$, $x(0)=0$ is a centered $H$-ss
continuous process and $\Delta$ is a star-shaped domain centered at $0$.
The following notation will be used:
$M_T=\sup \{ x(s), \, s\in T\cdot \Delta \}$ is the maximum of
$x$ in $T\cdot \Delta$;
$G_T=\{ s:\, x(s)=M_T, \, s\in T\cdot \Delta \}$ is the set of positions
of $M_T$;
$Z_T=\{ s:\, x(s)=0, \, s\in T\cdot \Delta \}$ is the
set of zeroes $x(s)$ in $s\in T\cdot \Delta$,
and $S_0=\{ s:\, E x(s)^2=0\}$ is the non-random set of zeroes  $x(s)$.
Note that $S_0=\{ 0\}$ for an~$H$-ss process with one-dimensional time.

The key observation for the analysis of (\ref{qq1}, \ref{qq2}) is that
a typical trajectory
of the $1$-$D$\, $H$-ss process $x(s)$ with small $M_T$ has its
$G_T$ and $Z_T$ close to $S_0$. For example, if $x(s)$ is
a Levy stable process and $\vert x\vert$ is not a subordinator,
then the position of maximum,
$G_T$, under the condition $\{ M_T<1\}$ has a non-trivial
limit distribution as $T\to \infty$ without any
normalization [8]. Therefore the events ${\cal A}_T$, namely,

\begin{eqnarray}\label{qq4}
\{ M_T<1\},\quad \{ G_T\subset U_0\},\quad
\{ Z_T\subset U_0,\quad x\vert_{\partial (T\cdot \Delta \setminus U_0)} <0 \},
\eqnum{4}
\end{eqnarray}
where $U_0$ is some vicinity of $S_0$
are considered as characteristic in the problem
(\ref{qq1}, \ref{qq2}).
By ${\cal A}_T(\xi)$, $\xi =M, G, Z$ we denote the above events
for further reference.
Consider a one-to-one mapping
$s_T : [0,1]^d\ni a\to s \in U_0\cap T\cdot \Delta$ and define

\begin{eqnarray*}
\delta^2_T(h) = \sup \{ E|x(s(a))-x(s(a'))|^2,\, a,a'\in [0,1]^d,\,
|a-a'|<h\}.
\end{eqnarray*}

{\bf Theorem~1.} Let $x(s)$, $x(0)=0$, $s\in R^d$ be a centered $H$-ss
continuous Gaussian process and
$\{ {\cal H}_x, \Vert  \cdot \Vert_x\}$ the reproducing kernel
Hilbert space associated with $x(\cdot)$ (see, e.g., [5]). Suppose there
exist a set $U_0 \supset S_0$, elements
$\varphi_T \in {\cal H}_x$,
and slowly increasing functions  $\psi (T)$ and  $\sigma_T$
such that

\begin{eqnarray}\label{qq5}
\nonumber
&\mbox {\rm (a)}&
\varphi_T(s) \ge 1,\quad  s\in (T\cdot \Delta ) \setminus U_0,\quad
\Vert  \varphi_T \Vert_x^2 = o(\psi (T)) \\
\nonumber
&\mbox {\rm (b)}&
\psi (T'a)/\psi (T)\to 1, \, \,T \to \infty \quad
\, \mbox {\rm for any}\,\,  a>0, \\
&\mbox {\rm \,\,}&
\noindent \mbox {\rm where}\, \,
T'=T/(\sigma^2_T\cdot \psi(T))^{1/2H}\quad \mbox {\rm and} \\
\nonumber
&\mbox {\rm (c)}&
\sigma_T \ge
\sup \{ (E |x(s)|^2)^{1/2},\,s\in U_0\cap T\cdot \Delta\}+
\bigg | \int\limits^1_0 \delta_T(h)\,d\sqrt{\ln 1/h}\bigg |.
\end{eqnarray}
Then all events ${\cal A}_T$ in (\ref{qq4}) have the decay exponents

\begin{eqnarray}\label{qq6}
\theta_\psi[x,\Delta] =
-\lim (\log\,P({\cal A}_T))/\psi (T), \quad T \to \infty
\eqnum{6}
\end{eqnarray}
simultaneously only and these exponents
are equal to each other.

{\it Remarks.} The condition (\ref{qq5}b) holds for
$\psi (T)$ and  $\sigma_T$ of the $(\log T)^\alpha$, $\alpha \ge 0$ type.
The conditions (\ref{qq5}b,c) are based on Fernique's estimate for the
distribution of the maximum of $|x(s)|$ in $U_0\cap T\cdot \Delta$.
If a similar bound due to Talagrand (see,e.g., [5]) is used, then
(5c) can be replaced with

\begin{eqnarray*}
\sigma_T \ge \sup \{ (E|x(s)|^2)^{1/2}, s\in U_0\cap T\cdot \Delta \} +
E\sup \{ |x(s)|, s\in U_0\cap T\cdot \Delta \}.
\end{eqnarray*}

P{\footnotesize ROOF}. For any
$\varphi_T \in \{ {\cal H}_x, \Vert  \cdot \Vert_x\}$
one has

\begin{eqnarray}\label{qq7}
\nonumber
{p}_T&=&
P\{ x(s)<1, \, \,  s\in T\cdot \Delta \} = \\
&=& P\{ x(s) - \varphi_T(s) < 1 - \varphi_T(s),\, \, s\in T\cdot \Delta \} =
E\,{\bf 1}_{\omega_T<1-\varphi_T}\pi (\omega_T),\qquad
\end{eqnarray}
where $\pi (\omega_T)$ is the Radon-Nikodym derivative of the two Gaussian
measures corresponding to the processes
$x(\cdot )+\varphi_T(\cdot )$ and $x(\cdot )$ on $T\cdot \Delta$.
By the Cameron-Martin formula, $\ln \pi (\omega_T)$ is a Gaussian
variable with mean \, $-\Vert \varphi_T \Vert_{x,T}^2/2$ \, \,and variance
$\Vert \varphi_T \Vert_{x,T}^2$, where

\begin{eqnarray*}
\Vert \varphi_T \Vert_{x,T}=
\inf \{ \Vert \widetilde{\varphi}_T \Vert_x,\,
\widetilde{\varphi}_T\in {\cal H}_x,
\widetilde{\varphi}_T(s)=\varphi_T (s), \,s\in T\cdot \Delta \} \le
\Vert \varphi_T \Vert_x\,\, .
\end{eqnarray*}

Therefore,
$E\,\pi^q(\omega_T) = \exp (q(q - 1)
\Vert  \varphi_{x,T} \Vert_x^2/2)$.

Applying to (7) H\"{o}lder's inequality with parameters

\noindent $(p,q)=((1-\varepsilon )^{-1}, \varepsilon^{-1})$,
we can continue (7):

\begin{eqnarray*}
{p}_T \le
P\{ x_T < 1-\varphi_T, \, \,  s\in T\cdot \Delta \}^{1-\varepsilon}
\exp (\Vert  \varphi_T \Vert_x^2 (\varepsilon^{-1}-1)/2).
\end{eqnarray*}

Taking into account (\ref{qq5}a), choose $\varepsilon$ as follows:
$\varepsilon^2=\Vert  \varphi_T \Vert^2/\psi(T)$.
Then $\varepsilon \to 0$ as $T \to \infty$ and

\begin{eqnarray}\label{qq8}
{p}_T\le  P\{ x(s)<0, \,
s\in T\cdot \Delta \setminus U_0\}^{1-\varepsilon}\,
e^{\varepsilon  \psi (T)/2} =
P\{ {\cal A}_T(Z)\}^{1-\varepsilon}\,e^{\varepsilon  \psi (T)/2}.
\eqnum{8}
\end{eqnarray}

The following inequality is obvious:

\begin{eqnarray}\label{qq9}
P({\cal A}_T(Z))\le  P ({\cal A}_T(G)).
\eqnum{9}
\end{eqnarray}

Lastly, for any $c_T>0$ one has

\begin{eqnarray*}
P({\cal A}_T(G)) &\le & P(G_T\subset U_0,\, M(U_0 \cap T\cdot \Delta ) < c_T) +
P(M(U_0 \cap T\cdot \Delta ) > c_T) \\
&\le & P(M_T<c_T) +  P(M(U_0\cap T\cdot \Delta) > c_T),
\end{eqnarray*}
where $M(I)=\max (|x(s)|,\,s\in I)$.

>From the $H$-ss property of $x(\cdot)$ and similarity of domains
$T\cdot \Delta$ we get

\begin{eqnarray*}
P(M_T<c_T) = P(M_{T'}<1), \quad T'=Tc_T^{-1/H}.
\end{eqnarray*}

By (\ref{qq5}b,c), for large enough
$T$ and $c_T=\sqrt{a\,\psi (T)}\sigma_T$ we have

\begin{eqnarray*}
P(M(U_0\cap T\cdot \Delta) > c_T) < \exp (-\bar{c}a\psi (T))
\end{eqnarray*}
with some constant $\bar{c}$ (see Fernique's estimate in [2]).
Therefore,

\begin{eqnarray}\label{qq10}
P({\cal A}_T(G)) \le P(M_{T'}<1) + e^{-ca\psi (T)}, \quad
T'=T(a \psi (T)\sigma_T^2)^{-1/2H}.
\eqnum{10}
\end{eqnarray}

We introduce the following notation:

\begin{eqnarray*}
\bar{\theta}_\xi =
\lim \sup \limits_{T \to \infty} -\log P({\cal A}_T(\xi))/\psi (T) \\
\underline{\theta}\,_\xi =
\lim \inf \limits_{T \to \infty} -\log P({\cal A}_T(\xi))/\psi (T)
\end{eqnarray*}
for  $\xi =M, Z, G$.

Taking into account (8) and (9), we obtain

\begin{eqnarray}\label{qq11}
\bar{\theta}_M \ge \bar{\theta}_Z \ge
\bar{\theta}_G \quad \mbox {\rm and} \quad
\underline{\theta}\,_M \ge \underline{\theta}\,_Z \ge
\underline{\theta}\,_G.
\eqnum{11}
\end{eqnarray}

Suppose that $\bar{\theta}_M<\infty$.
Using (10), we get

\begin{eqnarray*}
P(A_T(G)) \le P(M_{T'}<1)(1+\exp (-\bar{c}a\psi (T)+\theta \psi (T'))
\end{eqnarray*}
for $T>T_0$ and $\theta >\bar{\theta}_M$. By the assumption (\ref{qq5}b),
$\psi (T')=\psi (T)(1+o(1))$. Hence,

\begin{eqnarray}\label{qq12}
\bar{\theta}_G > \bar{\theta}_M \quad \mbox {\rm and} \quad
\underline{\theta}\,_G \ge \underline{\theta}\,_M
\eqnum{12}
\end{eqnarray}
because 'a' is arbitrary.

In the opposite case, \, $\bar{\theta}_M=\infty$,\,
$(\underline{\theta}\,_M=\infty)$
instead of (12) one has\,
$\bar{\theta}\,_G \ge ac$ \, $(\underline{\theta}\,_G \ge ac)$
for any~$a$.
Therefore, $\theta_\xi =\infty$ for $\xi =M, Z$ and $G$.
By (11), (12), it is clear that the decay exponents for
${\cal A}_T(\xi)$, $\xi =M,G,Z$ can exist only simultaneously and these
exponents are equal.
$\diamond$

Let $\Gamma =\{ g\}$ be a group of continuous $1:1$ transformations of $R^d$.
A random process $x(s)$ is homogeneous in $\Gamma$ if
$\{ x(gs)\} \stackrel{\rm d}{=} \{ x(s)\}$ for any $g\in \Gamma$.
Two domains $\Delta_1$ and $\Delta_2$ are equivalent if
$\Delta_2=g\Delta_1$ for some $g\in \Gamma$. By definition, a set of
equivalent domains $\{ \Delta_i \}$ form $\Delta_0$-parquet of $G$
if $\Delta_i \cap \Delta_j =\O$, $i\ne j$;
$\Delta_i=g_i\Delta_0$,  $g_i\in \Gamma$; and $[G]=[\cup \Delta_i]$, where
$[\Delta]$ is the closure of $\Delta$.

{\bf Proposition 1.} Let $x(s)$, $s\in R^d$ be a centered
continuous homogeneous in the group $\Gamma$  Gaussian process.
Suppose that
$Ex(t)x(s)>0$ on domains $\Delta_T \subset R^d$. If there exist integers
$r_1>r_2>1$ and $D$ such
that $\log r_1/\log r_2$ is irrational, for any $T>0$ and
$r=r_1,r_2$ the domain  $\Delta_T$ has $\Delta_{T/r}$-parquet of $r^D$
elements, then the limit

\begin{eqnarray*}
\theta := \lim_{T\to \infty}(-\log P(\Delta_T))T^{-D} < \infty
\end{eqnarray*}
exists. Here $P(\Delta_T) = P\{ x(s)<u, s\in \Delta_T \}$ and
$\Delta_T \subset \Delta_{T'}$ for $T<T'$.

P{\footnotesize ROOF}. If $\{ \Delta_i, i=1,...,r^D\}$ form the
$\Delta_{T/r}$-parquet of $\Delta_T$, then Slepian's lemma and
the homogeneity of $x$ in $\Gamma$ imply

\begin{eqnarray*}
P(\Delta_T) \ge \prod_i P(\Delta_i) = [P(\Delta_{T/r})]^q, \quad q=r^D
\end{eqnarray*}
or

\begin{eqnarray*}
\rho (T):= -\frac{\log P(\Delta_T)}{T^D}\le -\frac{\log P(\Delta_{T/r})}{(T/r)^D}.
\end{eqnarray*}
Therefore, for any $0<a\le r$, $r=r_1,r_2$ the sequence
$\rho (ar^n)$ tends from above to a limit $\theta_a(r)$ as
$n\to \infty$. We claim that $\theta_a(r_1)=\theta_1(r_2)$,
$0<a<r_1$. Indeed, the fractional parts of the numbers

\begin{eqnarray*}
\log_{r_2}a + n\rho, \quad n=0,1,...; \quad \rho = \log_{r_2}r_1
\end{eqnarray*}
are dense on $[0,1]$ because $\rho$ is irrational. Hence, there exist
$n_i^{\pm}$ and $m_i^{\pm}$ such that

\begin{eqnarray*}
T^+_i: = ar_1^{n_i^+} = r_2^{m_i^++\delta_i^+}, \quad
T^-_i: = ar_1^{n_i^-} = r_2^{m_i^--\delta_i^-},
\end{eqnarray*}
where $\delta_i^{\pm}\ge 0$ and $\delta_i^{\pm}\to 0$ as $i\to \infty$.
But for
$\widetilde{T}^{\pm}_i = r_2^{m_i^{\pm}}$

\begin{eqnarray*}
P(\Delta_{T^+_i}) \le P(\Delta_{\widetilde{T}^+_i}) \quad \mbox {\rm and}
\quad P(\Delta_{T^-_i}) \ge P(\Delta_{\widetilde{T}^-_i})
\end{eqnarray*}
i.e.,

\begin{eqnarray*}
\rho (T^+_i) \ge \rho (\widetilde{T}^+_i)
(\widetilde{T}^+_i/T^+_i)^D = \rho (\widetilde{T}^+_i)r_2^{-\delta^+_iD}.
\end{eqnarray*}
Therefore,

\begin{eqnarray*}
\theta_a(r_1) \ge \theta_1(r_2).
\end{eqnarray*}

Similarly, using $T^-_i$ and $\widetilde{T}^-_i$ we get

\begin{eqnarray*}
\theta_a(r_1) \le \theta_1(r_2).
\end{eqnarray*}

This concludes the proof.

\section{Fractional Brownian sheet}

{\bf Theorem~2.} 1. Let
$x(s)=x(s_1,s_2)$, $s_i>0$ be the fractional Brownian sheet
of index $0\le H<1$.

Consider $T\cdot \Delta_a =K_a\cap [0,T]^2$, where $K_a$ is the cone
$\{ s: 0<a<s_1/s_2 <a^{-1}\}$, $0\le a\le 1$. Then the limit

\begin{eqnarray}\label{qq13}
\lim_{T\to \infty} - \log\,P(x(s)<1, \, s\in T\cdot \Delta_a\}
/(\log\,T)^q = \theta_a(H)
\eqnum{13}
\end{eqnarray}
exists and is nontrivial $(0<\theta_a<\infty)$ for
$q =1+{\bf 1}_{a=0}$.

2. Consider the dual stationary sheet
$\xi (t)=\exp (-H(t_1+t_2))\,x(e^{t_1},e^{t_2})$.
Then

\begin{eqnarray}\label{qq14}
-\lim_{T\to \infty}\,T^{-2}\log P\{ \xi (t)<0, \,
t\in T\cdot \widetilde{\Delta}\} = \theta_0(H),
\eqnum{14}
\end{eqnarray}
where
$\widetilde{\Delta}=\{ t: t_1+t_2<2\} \cap [0,1]^2$.

{\it Remark.}
The relations (13) and (14) solve Problem 4 in ([4], p.227).

P{\footnotesize ROOF.}
By Li and Shao [4],

\begin{eqnarray*}
c_1(\log\,T)^2 < -\log\,P(x(s)<1, \, s\in [0,T]^2) < c_2(\log\,T)^2,
\end{eqnarray*}
where $c_i>0$. In fact, the proof of this relation leads to a more
general conclusion. Specifically, for the domain
$V_T=T\cdot \Delta_a \cap [1,T]^2$ we have

\begin{eqnarray*}
c_1n(T) < -\log\,P\{ x(s)<1, \, s\in V_T\} < c_2n(T),
\end{eqnarray*}
where $n(T)$ is the number of integer points in the domain
$\widetilde{V}_T=\{ (\log\,s_1,\,\log\,s_2) : (s_1,s_2)\in V_T\}$.
Obviously, $n(T)=\log\,T$ for the case $a>0$.
Hence we can expect nontrivial exponents with
$\psi (t)=\log T$ for $a>0$ and $\log^2T$ for $a=0$.

Consider $0<a\le 1$. To apply Theorem~1, we may use
$U_0=[0,1]^2\cap K_a$, $\psi (t)=\log T$ and

\begin{eqnarray*}
\varphi_T(s) = Ex(s_1,s_2)x(a,a)\cdot 4a^{-4H} =
\prod_{i=1}^2(1+(s_i/a)^{2H}-|s_i/a -1|^{2H}).
\end{eqnarray*}
One has $\varphi_T(s)\ge 1$ on $K_a\setminus U_0$ because
$\varphi_T(s)\ge 1$ on $[a,\infty )^2\supset K_a\setminus U_0$.
By the definition of $\varphi_T$,

\begin{eqnarray*}
\| \varphi_T\|^2 = E|x(a,a)\cdot 4a^{-4H}|^2 = 16a^{-4H} < \infty.
\end{eqnarray*}

It is clear that $E(x(s))^2<1$ on $U_0$.
Now $\delta_T^2(h)<ch^{2H}$, because

\begin{eqnarray*}
\delta^2(h)&=&E|x(s)-x(\widetilde{s})|^2 \le 2E|x(s)-x(\widehat{s})|^2
+ 2E|x(\widehat{s})-x(\widetilde{s})|^2 = \\
&=&2|s_1|^{2H}|s_2-\widetilde{s}_2|^{2H} +
2|\widetilde{s}_2|^{2H}|s_1-\widetilde{s}_1|^{2H} \le \\
&\le&2(|s_1|^{2H} + |\widetilde{s}_2|^{2H})\cdot |s_2-\widetilde{s}_2|^{2H},
\end{eqnarray*}
where $\widehat{s}=(s_1,\widetilde{s}_2)$.
Thus, all conditions in (5) are satisfied.

It remains to prove that the limit (6) exists for

\begin{eqnarray*}
{\cal A}_T = \{ x(s)<0, \,s\in T\cdot \Delta_a \setminus [0,1]^2\}
\end{eqnarray*}
with $\psi (T)=\log T$. In terms of the dual stationary process

\begin{eqnarray*}
\xi (t) = e^{-H(t_1+t_2)}x(e^{t_1}, e^{t_2})
\end{eqnarray*}
we have
${\cal A}_T = \{ \xi (t)<0, t\in \Delta_{\widetilde{T}_0, \widetilde{T}} \}$,
where $\widetilde{T}_0=0$, $\widetilde{T}=\ln T$ and

\begin{eqnarray*}
\Delta_{c,d} = \{ (t_1, t_2) : |t_1-t_2|<|\ln a|, \quad
t_1<d,\, \, t_2<d, \, \,(t_1\ge c \, \,
\mbox {\rm or} \, \,t_2\ge c)\}.
\end{eqnarray*}
We can apply Proposition~1 to the stationary process
$\xi (t)$, $t\in \Delta_{\widetilde{T}_0, \widetilde{T}}$
because for any integer $n$
$\Delta_{\widetilde{T}_0, \widetilde{T}}=\cup \Delta_{c_i, c_{i+1}}$,
$c_i=(\widetilde{T}/n)i$, $i=0,1,...,$ and
$\Delta_{c_i, c_{i+1}}=g(c_i)\Delta_{0, T/n}$, where
$g(a)(t_1,t_2)=(t_1+a, t_2+a)$. As a result, we get the existence
of the limit

\begin{eqnarray*}
\theta_a = -\lim (\log P({\cal A}_T))/\log T, \quad T\to \infty .
\end{eqnarray*}

Consider the case $a=0$. Here $T\cdot \Delta_0=[0,T]^2$ and
$S_0=\{ s: s_1s_2=0\}$. We will use $U_0=\{ s: s_1s_2<1\}$ as
a vicinity of $S_0$ and define the coordinates on
$U_0\cap [0,T]^2$ as follows:

\begin{eqnarray*}
a_1 = s_1/T, \quad a_2 = s_1s_2, \quad a=(a_1,a_2)\in [0,1]^2.
\end{eqnarray*}

One has $E|x(s)|^2\le 1$ on $U_0$. The next lemma shows that $\sigma_0^2(T)$
in (\ref{qq5}c) is $O(\log T)$ and therefore (\ref{qq5}b) holds.

{\bf Lemma~1.} For the fractional Brownian sheet

\begin{eqnarray*}
\sup \{ E|x(s)-x(\widetilde{s})|^2,\,
s,\widetilde{s} \subset U_0 \cap [0,T]^2, \,
|s_1-\widetilde{s}_1|<hT,\,
|s_1s_2-\widetilde{s}_1\widetilde{s}_2|<h\}  \\
= \delta_T^2(h) < c_H(1+2T^{-2}h^{-1})^{-2H}
\qquad \qquad \qquad \qquad
\end{eqnarray*}
and

\begin{eqnarray*}
-\int\limits_0^1\delta_T(h)\,d \sqrt{\ln h} \le c\sqrt{\ln T}, \quad
T>T_0.
\end{eqnarray*}

P{\footnotesize ROOF OF} L{\footnotesize EMMA}~1.
In the above estimate of $\delta^2(h)$ we may
use both $\widehat{s}=(s_1,\widetilde{s}_2)$ and
$\widehat{s}=(\widetilde{s}_1,s_2)$. Therefore,
taking into account the inequality

\begin{eqnarray*}
a^h + b^h \le c_h(a+b)^h, \quad c_h=2^{(1-h)_+}, \quad a>0,\,\, b>0, \,\,h>0,
\end{eqnarray*}
one has

\begin{eqnarray*}
E|x(s)-x(\widetilde{s})|^2 \le 2c_{2H}(E|w(s)-w(\widetilde{s})|^2)^{2H},
\end{eqnarray*}
where $w(s)$ is the Brownian sheet corresponding to the index $H=1/2$.
Consequently, it is sufficient to estimate $\delta_T^2(h)$ for $w(s)$ only.
One has

\begin{eqnarray*}
\delta^2 = E|w(s)-w(\widetilde{s})|^2 =
|(o\widetilde{s})\setminus (os)| + |(os)\setminus (o\widetilde{s})| = R_1+R_2,
\end{eqnarray*}
where $(os)=[0,s_1]\times [0,s_2]$ and $|G|$ is the area of $G$.

Consider $R_1=|o\widetilde{s}\setminus (os)|$ and suppose that
$s_1s_2=a_2>\widetilde{s}_1\widetilde{s}_2=\widetilde{a}_2$,
$a_2-\widetilde{a}_2=\delta a_2\le h$. There are three possibilities.

\noindent {\it First}, $(o\widetilde{s})\subset (os)$.
Then  $\delta^2=a_2-\widetilde{a}_2<h$.

\noindent {\it Second}, $\widetilde{s}_1>s_1$. In the notation
$s_1=Ta_1$, $\widetilde{s}_1=T\widetilde{a}_1$ one has

\begin{eqnarray*}
R_1 = (\widetilde{s}_1-s_1)\widetilde{s}_2 =
\frac{\widetilde{a}_2}{\widetilde{s}_1}  \delta s_1 =
\frac{a_2-\delta a_2}{s_1+\delta s_1}  \delta s_1 <
\frac{a_2 \cdot hT}{s_1+Th}
\end{eqnarray*}
because $\delta s_1=T\delta a_1<Th$. Since $s_1s_2=a_2$ and $s_2<T$,
we obtain $s_1>a_2/T$.
Hence $R_1 < h/(h+T^{-2})$ because $a_2<1$.

\noindent {\it Third}, $\widetilde{s}_2>s_2$,  $\widetilde{s}_1<s_1$.
Here $R_1=s_2(s_1-\widetilde{s}_1)=\delta s_1a_1/s_1$.

If $a_1/\widetilde{s}_1<T$, then
$a_1/(s_1-\delta s_1)<T$ or $s_1>a_1/T +\delta s_1$, so that

\begin{eqnarray*}
R_1 < \frac{a_1\delta s_1}{\delta s_1 + a_1/T} < h/(h+T^{-2}).
\end{eqnarray*}
If $a_1/\widetilde{s}_1>T$, then $\widetilde{a}_2/\widetilde{s}_1<1$,
i.e., $s_1>\delta s+\widetilde{a}_2/T$, and so

\begin{eqnarray*}
R_1 = \delta s_1\cdot a_1/s_1 <
\frac{\widetilde{a}_2+\delta a_2}
{\delta s_1+\widetilde{a}_2/T}  \delta s_1 <
\frac{h}{h+T^{-2}(1-h)}= \frac{hk}{h+T^{-2}k},
\end{eqnarray*}
where $k=(1-T^{-2})^{-1}<2$, $T>T_0$.
As a result, $R_1<2h(h+2T^{-2})^{-1}$.

Similarly, we get the same estimate for $R_2$. This completes the proof of
the first part of Lemma~1. Using for $x(\cdot)$ the estimate

\begin{eqnarray*}
\delta^2_T(h) \le c(1 + 2T^{-2}h^{-1})^{-2H},
\end{eqnarray*}
one has

\begin{eqnarray*}
-\int\limits_0^1\delta_T(h)\,d \sqrt{\ln 1/h} &\le&
c_1\int\limits_\varepsilon^\infty \sqrt{\ln u/\varepsilon}\,\,
(1+u)^{-H-1} du \\
&<&c_1\int\limits_0^\infty (\sqrt{|\ln u|} + \sqrt{\ln 1/\varepsilon}\,\, )
(1+u)^{-H-1}\,du = O(\ln 1/\varepsilon).
\end{eqnarray*}
We have integrated by parts and replaced $h$ by
$u=\varepsilon /h$ with $\varepsilon =2T^{-2}$.
$\diamond$

Let us continue the proof of Theorem~2. Consider the random variable

\begin{eqnarray*}
\eta_T = 4\sum_{|n|\le N}x(2^n, 2^{-n}), \quad N=[\log_2T]+1
\end{eqnarray*}
and the function $\varphi_T(s) = Ex(s)\eta_T$.

By definition, $\varphi_T(s) \in {\cal H}_x$ and
$\| \varphi_T\|^2_x=E\eta_T^2$.
The sequence
$\xi_n=4x(2^n, 2^{-n})$, $n=0,\pm 1,...,$ is stationary with the
correlation function

\begin{eqnarray*}
B(n) = (2^{-2Hn}+1 - |2^{-n}-1|^{2H})^2\cdot 2^{2Hn} > 0.
\end{eqnarray*}
Since $B(n)=O(2^{-n\gamma})$, $\gamma =2\cdot \min (H,1-H)$,
it follows that

\begin{eqnarray*}
\| \varphi_T\|^2_x = \sum_{(n,m)\le N}B(n-m) \le c\cdot N,
\end{eqnarray*}
where $c=\sum\limits_{-\infty}^\infty B(n)<\infty$.
Hence,  $\| \varphi_T\|^2_x=O(\ln T)$.

To show that $\varphi_T(s)>1$ on $[0,T]^2\setminus U_0$, note that

\begin{eqnarray*}
f(x) =|x|^{2H} + 1 - |x-1|^{2H}>1,\quad x>1/2, \quad H>0.
\end{eqnarray*}
But
$\psi_a(s) = Ex(s)\cdot 4x(a,a^{-1}) = f(s_1/a)f(s_2a)$.
Hence,\, $\psi_a(s)>1$ \,  for
$s\in [a/2, \infty)\times [a^{-1}/2, \infty)$ and

\begin{eqnarray*}
\varphi_T(s) = \sum_{|n|<N}\psi_{2^n}(s) > 1
\end{eqnarray*}
for $s\in \bigcup\limits_{|n|\le N} [2^{n-1}, \infty )
\times [2^{-n-1}, \infty )
\supset [0,T]^{2}\setminus U_0$ because $\psi_a(s)\ge 0$.

To prove (6), we consider

\begin{eqnarray*}
A_T(Z) = \{ x(s)<0, s\in [0,T]^2\setminus U_0\}
\end{eqnarray*}
and the dual stationary process $\xi (t)$ again. In terms of $\xi$

\begin{eqnarray*}
A_T(Z) = \{ \xi (t)<0, t_{1}<\widetilde{T},
t_{2}<\widetilde{T}, t_{1}+t_{2}>0\},\quad  \widetilde{T}=\ln T.
\end{eqnarray*}
Since

\begin{eqnarray}\label{qq15}
\{ \xi (e_1t_1+a_1, e_2t_2+a_2)\} \stackrel{\rm d}{=}\{ \xi (t_1,t_2)\}
\eqnum{15}
\end{eqnarray}
for any fixed $(a_1,a_2)$ and $e_i=\pm 1$, the event

\begin{eqnarray*}
A_{\widetilde{T}}(Z) = \{ \xi (t)<0, t_i\ge 0,
t_1+t_2<2\widetilde{T} \}
\end{eqnarray*}
has the same probability as $A_T(Z)$. To apply Proposition~1, we
divide

\begin{eqnarray*}
\Delta_{\widetilde{T}}= \{ (t_1, t_2):
t_1>0, t_2>0, t_1+t_2<2\widetilde{T}\}
\end{eqnarray*}
into $n^2$ equal triangles $\Delta_i$ using lines $t_1=a_ni$, $t_2=a_ni$
and $t_1+t_2=a_ni$, $i=1,2,...$ with $a_n=(2\widetilde{T})/n$.
By (15), all $\Delta_i$ are equivalent. Hence, we can apply
Proposition~1 to $\{ \xi (t), t\in \Delta_{\widetilde{T}}\}$ with
the parameter $D=2$.

\section{Integrated fractional Brownian motion}

Let $x_H(s)=\int\limits_0^sb_H(t)\,dt$, $s\in R^1$
be the integrated fractional Brownian motion and
$\Delta_i=(i,1)$, $i=0,-1$. Sinai [12] showed that

\begin{eqnarray*}
\lim_{T\to \infty} -\log P\{ x_H(s)<1, \,
s\in T\cdot \Delta_i\} /\log T = \theta_H(\Delta_i),
\end{eqnarray*}
exists for $H=1/2$ and $\theta_{1/2}(\Delta_0)=1/4$. The existence of both
$\theta_H(\Delta_0)$ for $0<H<1$ and $\theta_H(\Delta_{-1})$ for
$1/2\le H<1$ can be deduced from the relation
$Ex_H(s)x_H(t)\ge 0$ for $s,t\in \Delta_i$ (see [4] or
Proposition~1).

The numerical simulation [9, 11] suggests that
$\theta_0=H(1-H)$ and $\theta_{-1}=1-H$.
The symmetry of the hypothetical quantity
$\theta_0$, i.e., $\theta_0(H)=\theta_0(1-H)$, is quite unexpected.
For this reason the lower bound for $\theta_0$ given below is important.
As the numerical analysis shows, differences in the
behavior of $b_H(t)$ for $H<1/2$ and $H>1/2$ are seen in the
antisymmetry of the logarithmic correction
exponent, that is

\begin{eqnarray*}
P(x_H(t)<0, \, \, t\in [1,T]) \sim cT^{\theta_0}(\ln T)^{\alpha (H)},
\end{eqnarray*}
where $\alpha (H)=-\alpha (1-H)$ (for more detail see [11]).

{\bf Proposition~2.} For $x_H(s)$ there exists $0<\rho <1$ such that

\begin{eqnarray}\label{qq16}
\rho H(1-H) \le \theta_H(\Delta_0) \le \theta_H(\Delta_{-1}) \le 1-H
\eqnum{16}
\end{eqnarray}
if $\theta_H(\Delta_{-1})$ exists.

P{\footnotesize ROOF}. We begin by evaluating 
$\theta_{-1}:=\theta_H(\Delta_{-1})$ for the case $\Delta_{-1} =(-1,1)$.
Let $c(t)$ be the convex minorant of
$x(t)+t^2/2$. Then $\Dot{c}(t)$
is a continuous increasing function and $d \Dot{c}(t)$
is a measure with support $S$. If $t\in S$, then
$\Dot{c}(t)=\Dot{x}(t)$, where
$\Dot{x}=\mbox {\rm FBM}$; $t$ cannot be an
isolated point because of the continuity of $\Dot{c}(t)$.
If $\Delta \cap S \ne \O$, then

\begin{eqnarray*}
\int\limits_\Delta d\Dot{c}\,  =
|\Dot{c}(s_1) - \Dot{c}(s_2)| =
|\Dot{x}(s_1) - \Dot{x}(s_2) + s_1 - s_2|
\end{eqnarray*}
for some points $s_1,s_2\subset \Delta \cap S$. As is well known, FBM
belongs to the H\"{o}lder class of index $\alpha <H$, [3]. More exactly, for
any finite interval $[a,b]$ and $\varepsilon >0$\,
there exist $\delta_0>0$ and $c>0$ such that
$|\Dot{x}(t) - \Dot{x}(s)| < c|t - s|^{H-\varepsilon}$
a.s. as soon as  $|t-s|<\delta_0$. Hence,

\begin{eqnarray*}
\int\limits_\Delta d\Dot{c}\,  < b|\Delta |^{H-\varepsilon}
\end{eqnarray*}
for small enough $\Delta \subset [a,b]$ such that
$\Delta \cap S \ne \O$.
By Frostman's Lemma [3], $\mbox {\rm dim}S\ge H-\varepsilon$
a.s. Thus, $\mbox {\rm dim}S\ge H$ because $\varepsilon$
is arbitrary.

On the other hand, if $\theta_{-1}$ is the decay exponent for the event
${\cal A}_T(M, \Delta )$, $\Delta =(-1,1)$, then
$\mbox {\rm dim}S\le 1-\theta_{-1}$ (see [9]).
Hence, $H \le \mbox {\rm dim}S\le 1-\theta_{-1}$.
>From the relation
${\cal A}_T(M, (0,1)) \supset {\cal A}_T(M, (-1,1))$ one has
$\theta_0 \le \theta_{-1}$, where $\theta_0$ is related to $\Delta =(0,1)$.

To evaluate $\theta_0$, note that

\begin{eqnarray}\label{qq17}
\theta_0 = \theta_\xi =
-\lim_{T\to \infty} \ln P(\xi (t)<0, \,\, t\in (0,T))/T,
\eqnum{17}
\end{eqnarray}
where $\xi (t)=x(e^t)\cdot e^{-(1+H)t}$ is the dual stationary process.

\noindent By Lemma~2 (see below),

\begin{eqnarray*}
B_\xi(t) = E\xi (0) \xi (t) \le [\cosh (a(H)t)]^{-1},
\quad a(H)=\rho H(1-H)
\end{eqnarray*}
for some fixed $0<\rho <1$.

Applying Slepian's lemma [5], one has

\begin{eqnarray*}
P(\xi (t)<0, \,\, t\in (0,T)) \le  P(\eta (t)<0, \,\, t\in (0,2a(H)T)),
\end{eqnarray*}
where $\eta (t)$ is the stationary process with correlation function
$[\cosh(t/2)]^{-1}$. Hence  $\theta_0 =\theta_\xi \ge 2a(H)\theta_\eta$,
where $\theta_\eta$ is the decay exponent of type (\ref{qq17}) for the process
$\eta (t)$.
By [4], $2\theta_\eta \ge 0.2$. Thus,
$\theta_0 \ge 0.2\rho H(1-H)$.
$\diamond$

{\bf Lemma~2.} Let $\xi (t)$ be the dual stationary process of the
integrated fractional Brownian motion of index $H$. Then the correlation
function of $\xi (t)$

\begin{eqnarray}\label{qq18}
\nonumber
B_H(t) = [2(1+H)(e^{Ht}+e^{-Ht}) - (e^{(1+H)t}+e^{-(1+H)t}) + \\
(e^{t/2}-e^{-t/2})^{2H+2}]/(2+4H) \qquad  \qquad \qquad \qquad
\eqnum{18}
\end{eqnarray}
is nondecreasing and there exists $0<\rho <1$ such that

\begin{eqnarray}\label{qq19}
B_H(t) \le 1/\cosh (\rho H\bar{H}t), \quad \bar{H}=1-H
\eqnum{19}
\end{eqnarray}
for any $0<H<1$.

P{\footnotesize ROOF}.
Suppose that $B_H(t)$ is nondecreasing and (\ref{qq19}) holds for
$t \in [0,t_0] \cup [t_1, \infty )$ with a fixed $\rho =\rho_0$.
Then (19) holds for any $t$ with $\rho =\rho' =\rho_0t_0/t_1$.
Actually, for $t<t_1$  

\begin{eqnarray*}
B_H(t) \le B_H(t\cdot t_0/t_1) \le 1/\cosh (\rho' H\bar{H} t)
\end{eqnarray*}
and for $t>t_1$

\begin{eqnarray*}
B_H(t) \le 1/\cosh (\rho_0 H\bar{H} t) \le 1/\cosh (\rho' H\bar{H} t).
\end{eqnarray*}
Now we prove our admissions.

{\it Admission~1}: $B'_H(t) \le 0$. One has

\begin{eqnarray}\label{qq20}
B'_H(t) = -(1+H)(2+4H)^{-1} e^{(1+H)t} R(e^{-t}),
\eqnum{20}
\end{eqnarray}
where

\begin{eqnarray}\label{qq21}
R(x) = 1-2Hx-(1-x)^{2H}(1-x^2)+2Hx^{2H+1}-x^{2+2H}.
\eqnum{21}
\end{eqnarray}

Below we use the expansion

\begin{eqnarray}\label{qq22}
\nonumber
(1-x)^q &=& 1 + \sum_1^N q(q-1)...(q-k+1) (-x)^k/k! + \\
&+&q(q-1)...(q-N)(-x)^{N+1} r_{N+1}(q)/N!,
\eqnum{22}
\end{eqnarray}
where $(q\lor N)^{-1}\le r_N(q)\le (q\land N)^{-1}$.

Using (\ref{qq21}) and (\ref{qq22}), we get

\begin{eqnarray*}
R(x) &=& (1-2Hx)x^2-(2H-1)x^2(H+\bar{H}\cdot 2H r_3(2H)x)(1-x^2)+ \\
&+&x^{2H-1}(2H-x)x^2.
\end{eqnarray*}
In the case $2H>1$ one has $2Hr_3(2H)<1$ and $x^{2H-1}>x$. Hence,

\begin{eqnarray}\label{qq23}
R(x) \ge x^2(1-x^2)\bar{H}((2H-1)(1-x)+2)\ge 0.
\eqnum{23}
\end{eqnarray}
In the opposite case, $2H<1$, we have

\begin{eqnarray}\label{qq24}
\nonumber
R(x) &=& (1-x^{2H})x^2+2Hx^{2H+1}(1-x^{2\bar{H}})+ \\
&+&H(1-2H)x^2(1+2\bar{H}r_3(2H)x)(1-x^2)\ge 0.
\eqnum{24}
\end{eqnarray}
Thus, $R(x)\ge 0$ and $B_H(t)$ is nondecreasing.

{\it Admission~2}: relation (\ref{qq19}) for small $t$.

By (\ref{qq20}),

\begin{eqnarray}\label{qq25}
1 - B_H(t) \ge b(H) \int\limits_x^1 R(u)\, du, \quad x=e^{-t},
\eqnum{25}
\end{eqnarray}
where $b(H)=(1+H)/(2+4H)$. Using (\ref{qq23}) and (\ref{qq24}), one has

\begin{eqnarray}
\nonumber
R(u) \ge \begin{cases}
2\bar{H}u^2(1-u^2), \qquad \qquad \qquad \quad 2H>1 \\
(1-u^{2H})u^2 \ge 2Hu^2(1-u), \quad 2H < 1.
\end{cases}
\nonumber
\end{eqnarray}
Substituting these inequalities in (\ref{qq25}), we get on integration

\begin{eqnarray*}
1 - B_H(t) \ge b(1)\cdot 2\bar{H}(1-x)^2\varphi_1(x_0),
\quad 2H>1, \quad t<-\ln x_0=t_0,
\end{eqnarray*}
where $\varphi_1(x)=(2+4x+6x^2+3x^3)/15$, and

\begin{eqnarray*}
1 - B_H(t) \ge b(1/2)\cdot 2H(1-x)^2\varphi_2(x_0),
\quad 2H<1, \quad t<t_0,
\end{eqnarray*}
where $\varphi_2(x)=(1+2x+3x^2)/12$.

\noindent As a result, $1-B_H(t)\ge 2AH\bar{H}(1-x)^2$, $t<t_0$,
where

\begin{eqnarray*}
A = \min (\varphi_1(x_0)\cdot b(1),\,
\varphi_2(x_0)\cdot b(1/2)) = \varphi_1(x_0)/3.
\end{eqnarray*}

Since

\begin{eqnarray*}
1-x=1-e^{-t}>t\cdot e^{-t}>x_0t,
\end{eqnarray*}
$1-u^2<1/\cosh (\sqrt{2}\,u)$, and $1>2(H\bar{H})^{1/2}$, one has

\begin{eqnarray*}
B_H(t) \le 1/\cosh (2x_0(AH\bar{H})^{1/2}) < 1/\cosh (\rho H\bar{H}t),
\quad t<\ln 1/x_0=t_0,
\end{eqnarray*}
where $\rho =4x_0A^{1/2}$.

{\it Admission~3}: relation (\ref{qq19}) for large $t$.

Using (\ref{qq18}) and (\ref{qq22}), we get

\begin{eqnarray}\label{qq26}
\nonumber
B_H(t) \le (1+H)(1+2H)^{-1} e^{-Ht} - 2^{-1}(1+2H)^{-1} e^{-(1+H)t} +
\qquad \quad \\
(1+H) e^{-\bar{H}t}/2 - H(1+H) e^{-(1+\bar{H})t}/3
+ H(2H-1) e^{-t(2+\bar{H})}{\bf 1}_{2H>1}/6.\quad
\eqnum{26}
\end{eqnarray}

Here we use (\ref{qq22}) with $q=2H+2$, $N=1$ for the case $2H<1$ and
with $q=2H+2$, $N=2$ for the case $2H>1$.

{\it The case} $2H<1$. By (\ref{qq26}),

\begin{eqnarray}\label{qq27}
\nonumber
B_H(t) \le (1+H)(1+2H)^{-1} e^{-Ht} \times \qquad \\
\times [1 + 0.5 e^{-(1-2H)t}\{ 1+2H-(1+H)^{-1} e^{-2Ht}\} ].
\eqnum{27}
\end{eqnarray}
Using $1-e^x\le x$, we get

\begin{eqnarray*}
0,5\{ \cdot \} &=& 0.5\{ 2H+H(1+H)^{-1}+
(1+H)^{-1}(1-e^{-2Ht})\} \\
&\le & H(1+2^{-1}(1+H)^{-1}+t(1+H)^{-1})
\le H(1.5+t).
\end{eqnarray*}

If $2H<1-\varepsilon$, then

\begin{eqnarray*}
e^{-(1-2H)t}(1.5+t) \le
\max_t (e^{-\varepsilon t}(1.5+t)) = e^{1.5\varepsilon}/
(e\varepsilon) = a_\varepsilon.
\end{eqnarray*}

Hence, for $2H<1-\varepsilon$

\begin{eqnarray}\label{qq28}
B_H(t) \le e^{-Ht}(1+Ha_\varepsilon)<
e^{-H(t-a_\varepsilon)}.
\eqnum{28}
\end{eqnarray}
By (\ref{qq27}), for $1-\varepsilon <2H<1$

\begin{eqnarray}\label{qq29}
B_H(t) \le 2e^{-Ht} \le e^{-H(t-b_\varepsilon)},
\eqnum{29}
\end{eqnarray}
where $b_\varepsilon =\ln 4/(1-\varepsilon )$.

{\it The case} $2H>1$. By (\ref{qq26}),

\begin{eqnarray*}
B_H(t) &\le & 2^{-1}(1+H) e^{-\bar{H}t} + \\
&+&(1+H)(1+2H)^{-1} e^{-Ht}
(1-H(1+2H) e^{-2\bar{H}t}/3) - \\
&-&2^{-1}(1+2H)^{-1} e^{-(1+H)t}
(1- H(4H^2-1) e^{-2\bar{H}t}/3) \\
&\le &0.5(1+H) e^{-\bar{H}t}
[1+2(1+2H)^{-1} e^{-(2H-1)t} \times \\
&\times&\{ 1- H(1+2H)/3+H(1+2H)(1-e^{-2\bar{H}t})/3\} ].
\end{eqnarray*}
Arguing as above, we find that
$\{ \cdot \} <\bar{H}(1+2t)$ and $[\cdot ]<1+\bar{H}d_\varepsilon$
at last arriving at a bound like (\ref{qq28})
if $0<2\bar{H}<1-\varepsilon$. In the case $1-\varepsilon <2\bar{H}<1$,
we get (\ref{qq29}) again. Thus,

\begin{eqnarray*}
B_H(t) \le  e^{-\bar{H}\land H(t-t_0)}.
\end{eqnarray*}
But $(\bar{H}\land H)\cdot (t-t_0)>H\bar{H}t$ for
$t>t_0(1-H\lor \bar{H})^{-1}>2t_0$.

\noindent Hence, $B_H(t)\le 1/\cosh (H\bar{H}t)$ for $t>2t_0$.

\newpage

\bigskip

\noindent
I{\footnotesize NTERNATIONAL} I{\footnotesize NSTITUTE OF} \\
E{\footnotesize ARTHQUAKE} P{\footnotesize REDICTION}
T{\footnotesize HEORY AND} \\
M{\footnotesize ATHEMATICAL} G{\footnotesize EOPHYSICS}, \\
79, k2, Varshavskoe sh., 117556, Moscow, \\
RUSSIA \\
E-mail: molchan@mitp.ru

\end{document}